\pdfoutput=1
\RequirePackage{ifpdf}
\ifpdf 
\documentclass[pdftex]{sigma}
\else
\documentclass{sigma}
\fi

\numberwithin{equation}{section}

\newtheorem{Theorem}{Theorem}[section]
\newtheorem*{Theorem*}{Theorem}

\theoremstyle{definition}

\newtheorem{Remark}[Theorem]{Remark}

\usepackage{mathrsfs}

\begin{document}

\allowdisplaybreaks

\newcommand{\arXivNumber}{2511.21433}

\renewcommand{\PaperNumber}{060}

\FirstPageHeading

\ShortArticleName{Polynomials of the Askey Scheme as Clebsch--Gordan Coefficients}

\ArticleName{Polynomials of the Askey Scheme\\ as Clebsch--Gordan Coefficients}

\Author{Nicolas CRAMP\'E~$^{\rm a}$, Lo\"ic POULAIN D'ANDECY~$^{\rm b}$ and Luc VINET~$^{\rm c}$}

\AuthorNameForHeading{N.~Cramp\'e, L.~Poulain d'Andecy and L.~Vinet}

\Address{$^{\rm a)}$~CNRS - Universit\'e de Montr\'eal CRM-CNRS,\\
\hphantom{$^{\rm a)}$}~P.O. Box 6128, Centre-ville Station, Montr\'eal, H3C 3J7, Canada}
\EmailD{\mail{nicolas.crampe@cnrs.fr}}

\Address{$^{\rm b)}$~Laboratoire de math\'ematiques de Reims LMR, UMR 9008, \\
\hphantom{$^{\rm b)}$}~Universit\'e de Reims Champagne-Ardenne, Moulin de la Housse BP 1039,\\
\hphantom{$^{\rm b)}$}~51100 Reims, France}
\EmailD{\mail{loic.poulain-dandecy@univ-reims.fr}}

\Address{$^{\rm c)}$~IVADO and Centre de recherches math\'ematiques, Universit\'e de Montr\'eal,\\
\hphantom{$^{\rm b)}$}~P.O. Box 6128, Centre-ville Station, Montr\'eal, H3C 3J7, Canada}
\EmailD{\mail{luc.vinet@umontreal.ca}}

\ArticleDates{Received November 27, 2025, in final form June 12, 2026; Published online June 19, 2026}

\Abstract{Given a~semi-simple algebra equipped with a~coproduct, the Clebsch--Gordan coefficients are the elements of the transition matrices between direct product representation and its irreducible decomposition. It is well known that the Clebsch--Gordan coefficients of the Lie algebra $\mathfrak{sl}_2$ are given in terms of the dual Hahn polynomials. Taking the reversed point of view, we show that any finite dimensional family of polynomials belonging to the Askey scheme can be interpreted as Clebsch--Gordan coefficients of an algebra. The Hahn polynomials are thus associated to the oscillator algebra with the Krawtchouk polynomials treated through a~limit. The dual Hahn polynomials and Racah polynomials are seen to be associated to $\mathfrak{sl}_2$ with a~more general coproduct than the standard one. The $q$-Hahn polynomials are interpreted as Clebsch--Gordan coefficients of a~$q$-deformation of the oscillator algebra and the $q$-Racah polynomials are seen to be connected in this way to ${\rm U}_q(\mathfrak{sl}_2)$ with a~generalized coproduct.}


\Classification{33D45; 16T10}

\section{Introduction}

Clebsch--Gordan (CG) coefficients are fundamental in the representation theory of Lie algebras, where they enable the explicit direct sum decomposition of the tensor product of two irreducible representations. In physics, they arise as the coupling coefficients when two angular momenta are combined in quantum mechanics. The introduction of the coproduct, which is an algebra homomorphism from the algebra to the tensor product of the algebra, allows for the construction of tensor products of representations. This concept naturally leads to the search for non-cocommutative coproducts and, ultimately, to the notion of quantum groups~\cite{Drin1, Drin2}, which have since been intensively studied and applied in various contexts.

The structure of an algebra combined with its coproduct imposes strong constraints on its CG coefficients. For $\mathfrak{sl}_2$, this relationship yields recurrence relations that allows for identifying the CG coefficients as a~specific special function: the dual Hahn polynomials~\cite{Koor,Racah}.
The main goal of this paper is to perform the reverse procedure: starting with a~given set of CG coefficients, we derive the associated algebra $\mathcal{A}$ with the suitable homomorphism $\Delta\colon \mathcal{A}\to \mathcal{A}\otimes \mathcal{A}$.
Importantly, we do not require $\Delta$ to satisfy the full suite of axioms necessary for $\mathcal{A}$ to be a~Hopf algebra. In particular, the resulting homomorphisms are not necessarily coassociative. Consequently, $\Delta$ is referred to as a~generalized coproduct throughout this work.
 We prove that for any finite family of polynomials within the ($q$-)Askey scheme (see~\cite{KLS10}), there exists an algebra and a~generalized coproduct for which these polynomials serve as the corresponding CG coefficients.
Since dual Hahn polynomials belong to this scheme, our procedure not only recovers the $\mathfrak{sl}_2$ result but also yields a~more general coproduct for it. We also successfully identify the generalized coproducts for the ($q$-)oscillator algebra and the quantum group ${\rm U}_q(\mathfrak{sl}_2)$. The core technical tool enabling this procedure is the set of contiguity relations satisfied by the polynomials in the scheme, classified in~\cite{contiguity}.

{\bf Outline of the paper.} Section~\ref{sec:sl2} reviews the direct approach for $\mathfrak{sl}_2$: the well-known coproduct of $\mathfrak{sl}_2$ allows us to derive contiguity relations for the $\mathfrak{sl}_2$ Clebsch--Gordan coefficients, from which we deduce that these coefficients are proportional to the dual Hahn polynomials.
In Section~\ref{sec:GA}, the general reverse approach is described in detail, and we present the methodology to construct the algebra's defining relations and the generalized coproduct such that the associated CG coefficients are the polynomials we initially selected.
We then apply this procedure to the various finite families of the Askey scheme in Section~\ref{sec:askey} and of the $q$-Askey scheme in Section~\ref{sec:qaskey}.
Finally, Section~\ref{sec:outlook} is devoted to outlining the different future problems opened by the results obtained in this work.

{\bf Notations.} The generalized hypergeometric functions~\cite{GR,KLS10} are defined by
\[
{}_{r+1}F_r\left({{-n, a_1, a_2, \dots, a_{r} }\atop
{b_1, b_2, \dots, b_{r} }} \Bigg\vert  z\right)=\sum_{k=0}^{n}
\frac{(-n,a_1,a_2,\dots,a_r)_k}{k!(b_1,b_2,\dots,b_r)_k}z^k,
\]
for $r$, $n$ non-negative integers,
and where the Pochhammer symbols are
\[
(b_1,b_2,\dots,b_r)_k=(b_1)_k(b_2)_k\cdots (b_r)_k,\qquad (b)_k=b(b+1)\cdots (b+k-1).
\]
The $q$-hypergeometric functions~\cite{GR,KLS10} are defined by
\[
{}_{r+1}\phi_r\left({{q^{-n}, a_1, a_2, \dots, a_{r} }\atop
{b_1, b_2, \dots, b_{r} }} \Bigg\vert  q;z\right)=\sum_{k=0}^{n}
\frac{(q^{-n},a_1,a_2,\dots,a_r;q)_k}{(b_1,b_2,\dots,b_r,q;q)_k}z^k,
\]
for $r$, $n$ non-negative integers,
with the $q$-Pochhammer symbols given by
\[
(b_1,b_2,\dots,b_r;q)_k=(b_1;q)_k(b_2;q)_k\cdots (b_r;q)_k,\qquad (b;q)_k=(1-b)(1-qb)\cdots \bigl(1-q^{k-1}b\bigr).
\]
The $q$-binomial coefficients are
\[
 \left[\begin{array}{c}
 N \\
 n
\end{array}\right]_q =\frac{(q;q)_N}{(q;q)_n(q;q)_{N-n}}.
\]

\section[Clebsch--Gordan coefficients for sl\_2 revisited]{Clebsch--Gordan coefficients for $\boldsymbol{\mathfrak{sl}_2}$ revisited}\label{sec:sl2}

The Clebsch--Gordan coefficients for the Lie algebra $\mathfrak{sl}_2$ are proportional to the dual Hahn polynomials (see~\cite{Koor,NUV,Racah}). We show in this section how to prove this result using the contiguity and recurrence relations of these polynomials.

Let $E$, $F,$ and $H$ be the generators of the Lie algebra $\mathfrak{sl}_2$ satisfying the defining relations
$
 [H,E]=2E$, $ [H,F]=-2F$, $ [E,F]=H$.
The action of these elements on the lowest-weight Verma modules $V_\lambda$, for $\lambda\in\mathbb{C}$, with the basis $\{|\lambda,n\rangle | n=0,1,\dots\}$ is given by
\[H|\lambda,n\rangle=(\lambda+2n)|\lambda,n\rangle, \qquad E|\lambda,n\rangle=|\lambda,n+1\rangle, \qquad F|\lambda,n\rangle=-n(n+\lambda-1)|\lambda,n-1\rangle.\]
The representation $V_\lambda$ is irreducible if $\lambda$ is not a~negative integer: $\lambda\neq 0,-1,-2,\dots$.
The tensor product of two such representations $V_{\lambda_1}$ and $V_{\lambda_2}$ is isomorphic to a~direct sum of representations:%
\begin{equation}
 \label{eq:decompo}
V_{\lambda_1}\otimes V_{\lambda_2} \cong \bigoplus_{k\geq 0} V_{\lambda_1+\lambda_2+2k}.
\end{equation}
The parameters $\lambda_1$, $\lambda_2$ and $\lambda_1+\lambda_2$ are assumed not to be negative integers: $\lambda_1, \lambda_2, \lambda_1+\lambda_2\neq 0,-1,-2,\dots$, such that all the representations occurring in the previous formula are irreducible.

The action of the elements of the enveloping algebra $U(\mathfrak{sl}_2)$ on the modules $V_{\lambda_1}\otimes V_{\lambda_2}$ is given by the coproduct $\Delta$ which is an algebra homomorphism from $U(\mathfrak{sl}_2)$ to $U(\mathfrak{sl}_2)\otimes U(\mathfrak{sl}_2)$, given explicitly by
$\Delta(X)=X\otimes \mathbb{I}+\mathbb{I}\otimes X$, for $ X=E,F,H$,
with $\mathbb{I}$ the identity operator.

Knowing that $\Delta(H)=H\otimes \mathbb{I}+\mathbb{I}\otimes H$, a~vector $|\lambda_1+\lambda_2+2k,N-k\rangle$ belonging to the representation $V_{\lambda_1+\lambda_2+2k}$ in the direct sum of~\eqref{eq:decompo} can be expressed as the following linear combination
\begin{equation}\label{eq:defCG}
|\lambda_1+\lambda_2+2k,N-k\rangle=\sum_{n=0}^N \Gamma_n(k,N) |\lambda_1,n\rangle\otimes |\lambda_2,N-n\rangle, \ \end{equation}
with $N=0,1,\dots$ and $k=0,1,\dots,N$.
The above sum is taken over all the vectors in the eigenspace for $\Delta(H)$ with eigenvalue $\lambda_1+\lambda_2+2N$.
The coefficients $\Gamma_n(k,N)$ are called the Clebsch--Gordan coefficients.

{\bf Contiguity relations for CG coefficients.}
Acting with $\Delta(E)$ on relation~\eqref{eq:defCG} and projecting on $|\lambda_1,n\rangle\otimes |\lambda_2,N-n+1\rangle$, the following constraints for the CG coefficients are obtained, for $n=0,1,\dots N+1$,
\begin{equation}\label{eq:cont1}
 \Gamma_{n-1}(k,N)+\Gamma_n(k,N)=\Gamma_n(k,N+1) \qquad \text{with} \quad\Gamma_{-1}(k,N)=\Gamma_{N+1}(k,N)=0.
\end{equation}
Similarly, from the action of $\Delta(F)$ on~\eqref{eq:defCG} with $N\to N+1$, other constraints are obtained
\begin{equation}\label{eq:cont2}
 A_n \Gamma_{n+1}(k,N+1)+C_n \Gamma_n(k,N+1)=\mu(k)\Gamma_n(k,N),
\end{equation}
where
\begin{gather*}
A_n= (n+1)(n+\lambda_1),\qquad C_n=(N-n+1)(N-n+\lambda_2),\\ \mu(k)=(N-k+1)(N+k+\lambda_1+\lambda_2).
\end{gather*}
Both relations above are known as contiguity relations. Replacing in relation~\eqref{eq:cont2} the CG coefficients with argument $N+1$ using~\eqref{eq:cont1}, one obtains the following relation:
\[ A_n \Gamma_{n+1}(k,N)+(A_n+C_n)\Gamma_{n}(k,N)+C_n \Gamma_{n-1}(k,N)=\mu(k)\Gamma_n(k,N).\]
Considering carefully the coefficients, we can recognize $\Gamma_n(k,N)$ as being proportional to the dual Hahn polynomials (see, for example,~\cite{KLS10} for their recurrence relation)
\[
 \Gamma_n(k,N)=\binom{N}{n}\ {}_3F_2 \left({{-n, -k, k+\lambda_1+\lambda_2-1}\atop
{\lambda_1, -N}} \Bigg\vert  1\right).
\]
This hypergeometric function is indeed the dual Hahn polynomials $R_n(k(k+\lambda_1+\lambda_2-1);\allowbreak{\lambda_1-1},\lambda_2-1,N)$ in the usual notation (see, e.g., \cite{KLS10}).
We can check for consistency that the contiguity relations~\eqref{eq:cont1} and~\eqref{eq:cont2} are also satisfied (see~\cite{OV16} or the case called (dHII) in~\cite{contiguity}).

In conclusion, we have just seen that
the knowledge of the coproduct structure allows to compute the CG coefficients.
In this paper, we propose the reverse construction: we impose the Clebsch--Gordan decomposition with the Clebsch--Gordan coefficients given by any polynomials of the ($q$-)Askey scheme, and we ask whether there is an algebra with a~class of representations admitting a~tensor product structure realizing this Clebsch--Gordan decomposition. For the success of this procedure, the knowledge of the contiguity relations is crucial: their classification has been done in~\cite{contiguity}.
The next section explains in more details this procedure.
It turns out that surprisingly this construction works for all the finite sequences of polynomials of the ($q$-)Askey scheme: the different cases are presented in detail in the last sections of this paper.

\section{General approach }\label{sec:GA}

In this section, the elements $E$, $F$, $H$
generate an algebra $\mathcal{A}$ given with a~class of infinite-dimensional representations $V_\lambda$ indexed by a~complex number $\lambda$ (they will be lowest-weight Verma modules). The vectors of a~basis of $V_\lambda$ is denoted $\{|\lambda,n\rangle | n=0,1,\dots\}$ and we postulate the action of the generators as
\[H|\lambda,n\rangle=(\lambda+2n)|\lambda,n\rangle, \qquad E|\lambda,n\rangle=|\lambda,n+1\rangle, \qquad F|\lambda,n\rangle=\phi(\lambda,n)|\lambda,n-1\rangle.\]
The coefficients for the action of $E$ have been chosen equal to $1$ without loss of generality since we have the freedom in the normalization of the vectors $|\lambda,n\rangle$, whereas
the coefficients $\phi(\lambda,n)$ in the action of $F$ are to be determined and will define the algebra $\mathcal{A}$. Note that for now, the action of the generators ensure the relations
$[H,E]=2E$, $ [H,F]=-2F$.
Only the commutator of $E$ and $F$ remains to be determined. In all examples below, we will find that this commutator is a~constant (the oscillator algebra) or linear in $H$ (the Lie algebra $\mathfrak{sl}_2$), or a~$q$-deformations of these two situations.

Next, we consider the tensor product of two representations
$V_{\lambda_1}\otimes V_{\lambda_2}$
with the canonical product basis denoted $|\lambda_1,n\rangle\otimes |\lambda_2,m\rangle$, and we demand that
it decomposes as the following direct sum of representations
\[V_{\lambda_1}\otimes V_{\lambda_2}\cong\bigoplus_{k\geq 0} V_{\lambda_1+\lambda_2+2k}.\]
This decomposition is the same as that for the Lie algebra $\mathfrak{sl}_2$, as recalled in the previous section.
We also consider only the cases where the action of the Cartan element $H$ in the tensor product is the sum of the action in the two components, namely, $\Delta(H)=H\otimes \mathbb{I}+\mathbb{I}\otimes H$.
Therefore, for~$N$ a~non-negative integer and $k=0,1,\dots,N$, the vectors in the representation $V_{\lambda_1+\lambda_2+2k}$ are given by
\begin{equation}\label{CG-relation}
|\lambda_1+\lambda_2+2k,N-k\rangle=\sum_{n=0}^N P_n(k,N)|\lambda_1,n\rangle\otimes |\lambda_2,N-n\rangle,
\end{equation}
where the coefficients $\{P_n(k,N) | n=0,\dots,N\}$ will be the CG coefficients for the algebra $\mathcal{A}$.

Our goal in this paper is to find the algebra $\mathcal{A}$ with the associated homomorphism $\Delta\colon \mathcal{A}\to \mathcal{A}\otimes \mathcal{A}$, when $P_n(k,N)$ is a~finite sequence of polynomials (up to a~normalisation) of the ($q$-)Askey scheme.
 Namely, the possibilities are the Krawtchouk, the (dual) Hahn, the Racah polynomials together with their various $q$-versions.
The homomorphism $\Delta$ is defined by
\begin{gather*}
\Delta(H)|\lambda_1+\lambda_2+2k,N-k\rangle=(\lambda_1+\lambda_2+2N)|\lambda_1+\lambda_2+2k,N-k\rangle,\\
 \Delta(E)|\lambda_1+\lambda_2+2k,N-k\rangle=|\lambda_1+\lambda_2+2k,N-k+1\rangle,\\
 \Delta(F)|\lambda_1+\lambda_2+2k,N-k\rangle=\phi(\lambda_1+\lambda_2+2k,N-k)|\lambda_1+\lambda_2+2k,N-k-1\rangle.
\end{gather*}

\begin{Remark} \label{rem:iso}
As mentioned in the introduction, all the properties for $\Delta$, such that $\mathcal{A}$ becomes a~Hopf algebra, are not required and we referred to $\Delta$ as a~generalized coproduct. In particular, we do not require any coassociativity condition. Note nevertheless that the two ways of performing the tensor products
$(V_{\lambda_1}\otimes V_{\lambda_2})\otimes V_{\lambda_3}\cong V_{\lambda_1}\otimes(V_{\lambda_2}\otimes V_{\lambda_3})$
give isomorphic representations (from the explicit decompositions as direct sums). However, the identity map is not the isomorphism and this suggests non-trivial Drinfeld coassociators for our new tensor structures.
\end{Remark}

{\bf Strategy.} The idea is to use the so-called contiguity relations for the family of polynomials~$P_n(k,N)$. The contiguity relations which will be relevant for our purposes and in a~suitable notation are of the form
\begin{gather}\label{contiguity1}
P_n(k,N+1)=\alpha_1(n,N)P_{n-1}(k,N)+\alpha_2(n,N)P_n(k,N),\\
\label{contiguity2}
\mu_k(N)P_n(k,N-1)=\beta_1(n,N)P_{n+1}(k,N)+\beta_2(n,N) P_n(k,N).
\end{gather}
 Let us emphasize that the parameters (other than $N$ of course) in the polynomials remain constant in these relations. Such relations exist for all family of polynomials of the ($q$-)Askey scheme (see~\cite{contiguity} where they have been classified). The fact that these relations are closely related to our setting relies on the following general observations.

Given formula~\eqref{CG-relation} for the Clebsch--Gordan decomposition, one can calculate independently the action on these vectors of $E\otimes \mathbb{I}$, $\mathbb{I}\otimes E$ and $\Delta(E)$ as follows:
\begin{gather*}
\Delta(E)|\lambda_1+\lambda_2+2k,N-k\rangle = |\lambda_1+\lambda_2+2k,N+1-k\rangle\\
 \phantom{\Delta(E)|\lambda_1+\lambda_2+2k,N-k\rangle }{}= \sum_{n} P_n(k,N+1)|\lambda_1,n\rangle\otimes |\lambda_2,N+1-n\rangle,\\
(E\otimes \mathbb{I})|\lambda_1+\lambda_2+2k,N-k\rangle = \sum_{n} P_n(k,N)|\lambda_1,n+1\rangle\otimes |\lambda_2,N-n\rangle\\
 \phantom{(E\otimes \mathbb{I})|\lambda_1+\lambda_2+2k,N-k\rangle}{} = \sum_{n} P_{n-1}(k,N)|\lambda_1,n\rangle\otimes |\lambda_2,N+1-n\rangle,\\
(\mathbb{I}\otimes E)|\lambda_1+\lambda_2+2k,N-k\rangle = \sum_{n} P_n(k,N)|\lambda_1,n\rangle\otimes |\lambda_2,N+1-n\rangle.
\end{gather*}
The terms appearing in front of $|\lambda_1,n\rangle\otimes |\lambda_2,N+1-n\rangle$ are
$P_n(k,N+1)$, $ P_{n-1}(k,N)$, $ P_n(k,N)$,
and therefore a~relation among these three terms amount to a~relation giving $\Delta(E)$ in terms of~${E\otimes \mathbb{I}}$ and $E\otimes \mathbb{I}$. This is exactly what the first contiguity relation~\eqref{contiguity1} is doing. It turns out that in all the cases considered, we can express $\Delta(E)$ as
\[\Delta(E)=x(E\otimes \mathbb{I})+y(\mathbb{I}\otimes E),\]
where the coefficients $x$ and $y$ are determined to recover $\alpha_1(n,N)$ and $\alpha_2(n,N)$, respectively.
In Section~\ref{sec:Hahn}, these computations are performed in details for the Hahn polynomials.
These coefficients will often depend on the parameters of the polynomials we started with and also on the parameters $\lambda_1$, $\lambda_2$ of the representations. We will also formulate $x$ and $y$ as algebraic elements in $\mathcal{A}\otimes \mathcal{A}$ whose actions on $|\lambda_1,n\rangle\otimes |\lambda_2,N+1-n\rangle$ reproduce $\alpha_1(n,N)$ and $\alpha_2(n,N)$, respectively. This will involve the Cartan elements $H\otimes \mathbb{I}$ and $\mathbb{I}\otimes H$ (as they do, for example, for~${\rm U}_q(\mathfrak{sl}_2)$) and also a~central Casimir element $C$ of the algebra $\mathcal{A}$. In order for this procedure to work, it will be sometimes necessary to divide by algebraic expressions so strictly speaking we will be working in some localization of $\mathcal{A}\otimes \mathcal{A}$ (see the remarks below for more details).

Similar reasoning is also valid for $\Delta(F)$, $F\otimes \mathbb{I}$ and $\mathbb{I} \otimes F$. The only difference is that the action of $F$ involves the function $\phi$ we look for and it must be given by
\begin{equation}\label{eq:phimu}
 \phi(\lambda_1+\lambda_2+2k,N-k)=\mu_k(N),
\end{equation}
where $\mu_k(N)$ is the function in the left-hand side of the contiguity relation~\eqref{contiguity2}. Let us emphasize that a~proportionality factor could be added in the previous relation. This factor would change the commutation relation but the resulting algebra is isomorphic to the one with this factor put to one by a~rescaling of the generator $F$.

The coefficients $x'$ and $y'$ appearing in
$\Delta(F)=x'(F\otimes \mathbb{I})+y'(\mathbb{I}\otimes F)$,
are determined so that their action on $|\lambda_1,n\rangle\otimes |\lambda_2,N-1-n\rangle$ provides \smash{$\frac{\beta_1(n,N)}{\phi(\lambda_1,n+1)}$} and \smash{$\frac{\beta_2(n,N)}{\phi(\lambda_2,N-n)}$}, respectively.

We emphasize that, having prescribed the actions of $\Delta(H)$, $\Delta(E)$ and $\Delta(F)$ on the vectors~${|\lambda_1+\lambda_2+2k,N-k\rangle}$, it is a~consequence of our construction that these elements satisfy the relations of the algebra $\mathcal{A}$. In other words, it is automatic that the map $\Delta$ is an algebra homomorphism.

To summarize, this procedure allows to determine the algebra $\mathcal{A}$ with the determination of~$\phi$ and provide the associated homomorphism $\Delta$ with the computations of $x$, $y$, $x'$ and $y'$.

{\bf Orthogonality.} Let us choose a~scalar product in $V_{\lambda}$ such that
$
 \langle \lambda,m|\lambda,n \rangle=\delta_{n,m}$,
and consider the usual scalar product in $V_{\lambda_1}\otimes V_{\lambda_2}$ as the product of the scalar product in each space. We renormalize this scalar product by setting their norm $||\cdot ||$ to
\[\big|\big| |\lambda_1,n\rangle\otimes |\lambda_2,N-n\rangle\big|\big|^2=\Omega_{n,N},\]
and we choose $\Omega_{n,N}$ such that
\[\sum_{n=0}^NP_n(k,N)P_n(\ell,N)\Omega_{n,N}=\delta_{k,l}\Omega'_{\ell,N} .\]
For a~family of polynomials of the ($q$-)Askey scheme, the values of $\Omega_{n,N}$ and $\Omega'_{\ell,N}$ are known~\cite{KLS10}. With this scalar product, the vectors $|\lambda_1+\lambda_2+2\ell,N-\ell \rangle$ defined by~\eqref{CG-relation} are orthogonal and can be renormalized thanks to $\Omega'_{\ell,N}$. Note that the factor $\Omega_{n,N}$ does not always factorize as a~product of a~function of $n$ and a~function of $N-n$ and therefore cannot always be absorbed in a~renormalization of the vectors $|\lambda,n\rangle$.

{\bf Remarks on the generalized coproducts.} As mentioned before, the generalized coproducts $\Delta$ we will obtain is not always well defined from $\mathcal{A}$ to $\mathcal{A}\otimes \mathcal{A}$. The reason is that we need to divide by some elements, i.e., for example, by some polynomials in $H\otimes \mathbb{I}$ or in a~central element $\mathbb{I}\otimes C$. Therefore a~precise ``universal'' definition of $\Delta$ should involve some localization or some extension of $\mathcal{A}\otimes\mathcal{A}$. Another point of view would be to define, instead of this universal generalized coproduct, the following algebra homomorphism
\[\Delta_{\lambda_1,\lambda_2}\colon \ \mathcal{A} \to \operatorname{End}(V_{\lambda_1}\otimes V_{\lambda_2}),\]
depending on the parameters $\lambda_1$, $\lambda_2$. This is indeed sufficient for having well-defined tensor products of representations and a~well-posed CG problem. This homomorphism is valid for generic parameters, ensuring that the denominators appearing in its definition are invertible.

{\bf Remarks on the choice of the normalization.} The validity of the procedure presented above does not depend on the normalization of the CG coefficients we start with, but the generalized coproduct obtained can have a~different expression. Indeed, if we consider in relation~\eqref{CG-relation} the CG coefficient $f_k(N)g_n(N)P_n(k,N)$ instead of $P_n(k,N)$. The normalization factor~$g_n(N)$ could modify the coefficients in the contiguity relations $\alpha_1(n,N)$, $\alpha_2(n,N)$, $\beta_1(n,N)$ and~$\beta_2(n,N)$ and the obtained coproduct could have a~different expression. In the following, the normalization factor $g_n(N)$ is chosen such that the necessary localizations, discussed in the previous remark, are minimal. The normalization factor $f_k(N)$ could change the factors in the left-hand side of the contiguity relations as $\mu_k(N)$. It is chosen such that this factor is $1$ in the first relation.
There is also the freedom to multiply the second contiguity relation by a~constant. It would change the factor $\mu_k(N)$ and the obtained commutation relations. However, this simply corresponds to a~different normalization of $F$ and produces an isomorphic algebra.

\section{Askey scheme} \label{sec:askey}

In this section, we provide the explicit examples associated to the finite polynomial sequences of the Askey scheme. We start with the Hahn polynomials and show that it is associated to the oscillator algebra. We also perform a~limit to obtain the Krawtchouk polynomials. Then, we discuss the dual Hahn polynomials and obtain a~more general coproduct than the usual one discussed in Section~\ref{sec:sl2}. We end this section by applying the method to the Racah polynomials.

\subsection{Hahn polynomials and oscillator algebra} \label{sec:Hahn}

Let us consider for the functions appearing in~\eqref{CG-relation} the polynomials defined by
\[Q_n(k,N)=\binom{N}{n}\ {}_3F_2\left(\begin{array}{c} -n, n+\alpha+\beta-N+1, -k \\ \alpha+1, -N\end{array}\Bigg\vert 1\right),\]
which are proportional to the Hahn polynomials.
Their contiguity relations called (HIV) in~\cite{contiguity} are\footnote{We have replaced $\beta$ by $\beta-N$ compared to~\cite{contiguity}.}
\begin{gather*}
Q_n(k,N+1)=\frac{n+\alpha+\beta+1}{2n+\alpha+\beta-N}Q_{n-1}(k,N)+\frac{n+\alpha+\beta-N}{2n+\alpha+\beta-N}Q_n(k,N),\\
(N-k)Q_n(k,N-1)=\frac{(n+1+\alpha)(n+1)}{2n+2+\alpha+\beta-N}Q_{n+1}(k,N)\\
\phantom{(N-k)Q_n(k,N-1)=}{}+\frac{(n+1+\beta-N)(N-n)}{2n+2+\alpha+\beta-N}Q_n(k,N).
\end{gather*}

{\bf Oscillator algebra.} Looking at the left-hand side of the last equation, we take
${\phi(\lambda,k)=-k}$.
We have multiplied the second contiguity relation by a~global factor $-1$ for convenience (this does not change the algebra up to a~rescaling of the generator $F$). Therefore, the representations of $\mathcal{A}$ are
$H|\lambda,n\rangle=(\lambda+2n)|\lambda,n\rangle$, $E|\lambda,n\rangle=|\lambda,n+1\rangle$, $F|\lambda,n\rangle=-n|\lambda,n-1\rangle$,
and the defining relations of the algebra are
\begin{equation}\label{eq:osc}
[H,E]=2E, \qquad [H,F]=-2F, \qquad [E,F]=1.
\end{equation}
Therefore, the algebra $\mathcal{A}$ is the oscillator algebra $osc$. The Casimir element of this algebra is
$C=2EF+H$, with $ C|\lambda,n\rangle=\lambda|\lambda,n\rangle $.

{\bf Generalized coproduct.} As explained in the previous section to determine the generalized coproduct, we look for an element $x$ of the algebra $\mathcal{A}\otimes\mathcal{A}$ such that
\begin{equation}\label{eq:xH}
 x|\lambda_1,n\rangle\otimes |\lambda_2,N+1-n\rangle=\frac{n+\alpha+\beta+1}{2n+\alpha+\beta-N}|\lambda_1,n\rangle\otimes |\lambda_2,N+1-n\rangle.
\end{equation}
Knowing the actions of $H^{(1)}=H\otimes \mathbb{I}$, $H^{(2)}=\mathbb{I} \otimes H$, $C^{(1)}=C\otimes \mathbb{I}$ and $C^{(2)}=\mathbb{I} \otimes C$ on this vector, we can show that the following combination
\[
 x=\frac{1}{D}\bigl(H^{(1)}-C^{(1)}+2\alpha+2\beta+2\bigr),\]
where \smash{$D=H^{(1)}-H^{(2)}-C^{(1)}+C^{(2)}+2\alpha+2\beta+2$}, satisfies~\eqref{eq:xH}. As mentioned previously, we must consider the localization of $\mathcal{A}\otimes\mathcal{A}$ where $D$ is invertible for $x$ to be well defined.

The same computation for $y$ leads to the following formula for the generalized coproduct
\begin{equation}\label{eq:DEdual}
 \Delta(E)=\frac{1}{D}\bigl(\bigl(H^{(1)} -C^{(1)}+2\alpha+2\beta+2\bigr)(E\otimes \mathbb{I})+\bigl(C^{(2)} -H^{(2)}+2\alpha+2\beta+2\bigr)(\mathbb{I}\otimes E)\bigr).
\end{equation}
Similarly, we compute
\begin{equation}\label{eq:DFdual}
\Delta(F)=\frac{1}{D}\bigl(\bigl(H^{(1)}-C^{(1)}+2\alpha+2\bigr)(F\otimes \mathbb{I})+\bigl(C^{(2)}-H^{(2)}+2\beta\bigr)(\mathbb{I}\otimes F)\bigr).\end{equation}
 Replacing $C^{(i)}$ by $\lambda_i$, we can also see the formulas above as providing a~family of homomorphisms
\[\Delta_{\lambda_1,\lambda_2}\colon\ \mathcal{A} \to \operatorname{End}(V_{\lambda_1}\otimes V_{\lambda_2}) .\]
Let us emphasize that the construction depends on the two parameters $\alpha$ and $\beta$ which can be functions of $C^{(1)}$ and $C^{(2)}$ (or of $\lambda_1$ and $\lambda_2$). This is well defined under the restriction $\alpha+\beta\notin\mathbb{Z}$, because of the denominators.

Note that an automatic consequence of our construction is that the generalized coproduct, defined by relations~\eqref{eq:DEdual} and~\eqref{eq:DFdual}, is a~homomorphism from $osc$ to the correct localization of~$osc \otimes osc$, or to the endomorphism algebra $\operatorname{End}(V_{\lambda_1}\otimes V_{\lambda_2})$. Nevertheless, this can also be checked by a~direct calculation.

{\bf Limit to Krawtchouk polynomials.} Take $\alpha=p z$ and $\beta=(1-p)z$, and consider the limit $z\to+\infty$ in the previous computations. The polynomials $Q_n(k,N)$ become proportional to the Krawtchouk polynomials, namely they become
\[K_n(k,N)=\binom{N}{n}\ {}_2F_1\left(\begin{array}{c} -n, -k \\ -N\end{array}\Bigg\vert p^{-1}\right) .\]
The algebra remains the oscillator algebra defined by~\eqref{eq:osc}. However, the associated homomorphism simplifies and becomes
$
\Delta(E)=E\otimes \mathbb{I}+\mathbb{I} \otimes E$, $
\Delta(F)=pF\otimes \mathbb{I}+(1-p)\mathbb{I}\otimes F$.
Note that for this generalized coproduct no localization is needed.
\begin{Remark}\label{rema_coassociativitypq}
Let us denote by $\Delta_p$ the map above.
This generalized coproduct is not coassociative in general
\begin{equation}\label{eq:Dp}
 (\Delta_q\otimes Id)\circ\Delta_p \neq (Id\otimes \Delta_{q'})\circ\Delta_{p'}.
\end{equation}
However, as already mentioned in Remark~\ref{rem:iso}, they lead to isomorphic representations.
It is worth to notice that if the parameters in~\eqref{eq:Dp} satisfy
\begin{equation}\label{coassociativitypq}
p'=pq \qquad \text{and}\qquad 1-p=(1-p')(1-q'),\end{equation}
we have equality in~\eqref{eq:Dp}.
This suggests that the non-trivial Drinfeld coassociator reduces to the identity when conditions~\eqref{coassociativitypq} is satisfied.
\end{Remark}

{\bf Comparison with standard results~\cite{CVV,Zhosc}.} To associate a~coproduct with the oscillator algebra, one considers usually a~centrally extended version of it
\begin{equation}\label{eq:oscext}
[H,E]=2E,\qquad [H,F]=-2F,\qquad [E,F]=c,
\end{equation}
with $c$ a~central element.
For this latter algebra, the map $\Delta(X)=X\otimes \mathbb{I} + \mathbb{I} \otimes X$ (with $X=E,F,H,c$) is a~coproduct and the CG coefficients are also the Krawtchouk polynomials, with parameters $p$ given in terms of the values of the central elements $c_1=c\otimes \mathbb{I}$ and~${c_2=\mathbb{I} \otimes c}$ by
$ p=\frac{c_1}{c_1+c_2}$.
Furthermore, this coproduct is obviously coassociative.
Although the Krawtchouk polynomials appear in this
standard setting, the limit procedure presented in the above paragraph produces a~generalized coproduct $\Delta_p$ exhibiting a~broader class of tensor-product structures and without the necessity to centrally extend the oscillator algebra.

In this case, the four parameters $p$, $p'$, $q$, $q'$ appearing in the different Krawtchouk polynomials, when one considers the two ways of recoupling a~three-fold tensor product, are given in terms of the values of the central elements $c_1$, $c_2$ and $c_3=\mathbb{I}\otimes\mathbb{I}\otimes c$ by
\[q=\frac{c_1}{c_1+c_2}, \qquad p=\frac{c_1+c_2}{c_1+c_2+c_3}, \qquad q'=\frac{c_2}{c_2+c_3}, \qquad p'=\frac{c_1}{c_1+c_2+c_3}.\]
They are not independent, and in fact they satisfy exactly the constraints~\eqref{coassociativitypq}. This shows that with the homomorphism $\Delta_p$, we obtain a~recoupling depending on 4 independent parameters at the price of relaxing the coassociativity of the tensor product.

\subsection{Dual Hahn polynomials revisited}

In Section~\ref{sec:sl2}, we recall that the dual Hahn polynomials are the CG coefficients for the Lie algebra~$\mathfrak{sl}_2$.
Applying our procedure, we show that not only do we recover this result but we generalize it: there is one parameter which remains free in the CG coefficient.

Let us consider for the functions appearing in~\eqref{CG-relation} the polynomials defined by
\[
R^{\rm (dH)}_n(k,N)=\binom{N}{n}\ {}_3F_2\left(\begin{array}{c} -k, k+\alpha+\beta-N+1, -n \\ \alpha+1, -N\end{array}\Bigg\vert 1\right)
\]
with their contiguity relations called (dHII) in~\cite{contiguity}
\begin{gather*}
R^{\rm (dH)}_n(k,N+1)=R^{\rm (dH)}_{n-1}(k,N)+R^{\rm (dH)}_n(k,N),\\
(k-N)(N+k+\alpha+\beta+1)R^{\rm (dH)}_n(k,N-1)\nonumber\\
\qquad=-(n+1+\alpha)(n+1)R^{\rm (dH)}_{n+1}(k,N)+(N-n+\beta)(n-N)R^{\rm (dH)}_n(k,N).
\end{gather*}
The left-hand side of the second contiguity relation imposes
\[\phi(\lambda_1+\lambda_2+2k,N-k)=(k-N)(N+k+\alpha+\beta+1),\]
which implies that
$\phi(\lambda,k)=-k(k+\lambda+s)$,
for a~parameter $s$ with the constraint
$\alpha+\beta+1=\lambda_1+\lambda_2+s$.
In this case, the commutation relation between $E$ and $F$ is
$ [E,F]=H+1+s$.
By setting $s=-1$, this defining relation becomes $[E,F]=H$ and the algebra $\mathcal{A}$ becomes the Lie algebra $\mathfrak{sl}_2$.

The associated homomorphism $\Delta_{\lambda_1,\lambda_2}\colon\mathcal{A} \to \operatorname{End}(V_{\lambda_1}\otimes V_{\lambda_2})$ is given by
\begin{gather*}
 \Delta_{\lambda_1,\lambda_2} \colon\ E\mapsto E\otimes \mathbb{I}+\mathbb{I}\otimes E,\\
 F\mapsto \frac{H^{(1)}-\lambda_1+2\alpha+2}{H^{(1)}+\lambda_1+2s+2}F\otimes \mathbb{I}+\frac{H^{(2)}-\lambda_2+2\beta+2}{H^{(2)}+\lambda_2+2s+2} \mathbb{I}\otimes F.
\end{gather*}
In the previous expression, we have used the same notation for the abstract generators of the algebra and their representation.
As in Section~\ref{sec:sl2}, we assume that $\lambda_i+s+1\neq0,-1,-2,\dots $ such that the denominators in the right-hand side of the previous map do not vanish. Again the fact that the map $\Delta_{\lambda_1,\lambda_2}$ is an algebra homomorphism is a~consequence of our construction but can also be checked by a~direct calculation.

To pass from the map $\Delta_{\lambda_1,\lambda_2}$ to $\Delta$, a~generalized coproduct $\mathcal{A}\to \mathcal{A}\otimes\mathcal{A}$, we need to express~$\lambda_1$,~$\lambda_2$ in terms of the Casimir element $C$ of $\mathfrak{sl}_2$, which is given by
\begin{equation}\label{eq:Cassl2}
 C=4EF+H^2-2H \qquad \text{with} \quad C|\lambda,n\rangle=\lambda(\lambda-2)|\lambda,n\rangle .
\end{equation}
This can be done by extending $\mathcal{A}$ with a~central element $X$ satisfying $X(X-2)=C$. This generalized coproduct is well defined if one considers the localization where $H^{(i)}+X^{(i)}+2s+2$ are invertible.

Finally, note that we recover the usual coproduct of $\mathfrak{sl}_2$ when $s=-1$, $\alpha+1=\lambda_1$ and $\beta+1=\lambda_2$.

\subsection[Racah polynomials and generalized coproduct for sl\_2]{Racah polynomials and generalized coproduct for $\boldsymbol{\mathfrak{sl}_2}$}\label{ssec:Racah}

Let us consider for the functions appearing in~\eqref{CG-relation} the polynomials defined by\footnote{We have replaced $\beta$ by $\beta-N$ and $\gamma$ by $\gamma+N$ compared to~\cite{contiguity}.}
\[R^{\rm (R)}_n(k,N)=\binom{N}{n}{}_4F_3\left(\begin{array}{c} -n, n+\alpha+\beta-N+1, -k, k+\gamma \\ \alpha+1, \beta+\gamma+1, -N\end{array}\Bigg\vert 1\right),\]
which are proportional to the Racah polynomials.

The contiguity relations (see (RI) in~\cite{contiguity}) are
\[
R^{\rm (R)}_n(k,N+1)=\frac{n+\alpha+\beta+1}{2n+\alpha+\beta-N}R^{\rm (R)}_{n-1}(k,N)+\frac{n+\alpha+\beta-N}{2n+\alpha+\beta-N}R^{\rm (R)}_n(k,N),
\]
and
\begin{gather*}
(k-N)(N+k+\gamma)R^{\rm (R)}_n(k,N-1)\\
\qquad=-\frac{(n+\beta+\gamma+1)(n+\alpha+1)(n+1)}{2n+2+\alpha+\beta-N}R^{\rm (R)}_{n+1}(k,N)\\
\phantom{\qquad=}{}+\frac{(n+1+\beta-N)(n+1+\alpha-\gamma-N)(N-n)}{2n+2+\alpha+\beta-N}R^{\rm (R)}_n(k,N).
\end{gather*}
The left-hand side of the second contiguity relation imposes
\[\phi(\lambda_1+\lambda_2+2k,N-k)=(k-N)(N+k+\gamma),\]
which implies that
$\phi(\lambda,k)=-k(k+\lambda+s)$,
 with the constraint
$\gamma=\lambda_1+\lambda_2+s$.
As for the dual Hahn case, by choosing $s=-1$, the commutation relation between $E$ and $F$ is
$[E,F]=H$,
and we obtain the Lie algebra $\mathfrak{sl}_2$.

The resulting homomorphism $\Delta_{\lambda_1,\lambda_2}\colon \mathcal{A} \to \operatorname{End}(V_{\lambda_1}\otimes V_{\lambda_2})$ is given by
\begin{gather*}
\Delta_{\lambda_1,\lambda_2}\colon\ E\mapsto \frac{1}{D}\bigl(\bigl(H^{(1)}-\lambda_1+2\alpha+2\beta+2\bigr)E\otimes \mathbb{I}+\mathbb{I}\otimes \bigl(\lambda_2-H^{(2)}+2\alpha+2\beta+2\bigr)E\bigr),\\
 F\mapsto \frac{1}{D}\left(\frac{\bigl(H^{(1)}-\lambda_1+2\alpha+2\bigr)\bigl(H^{(1)}-\lambda_1+2\beta+2\gamma+2\bigr)}{H^{(1)}+\lambda_1+2s+2}F\otimes \mathbb{I}\right. \\
 \phantom{ F\mapsto }{}\left. +\frac{\bigl(\lambda_2-H^{(2)}+2\beta\bigr)\bigl(H^{(2)}-\lambda_2-2\alpha+2\gamma\bigr)}{H^{(2)}+\lambda_2+2s+2}\mathbb{I}\otimes F\right),\end{gather*}
with $D=H^{(1)}-H^{(2)}-\lambda_1+\lambda_2+2\alpha+2\beta+2$. In the previous map, one keeps the same notations for the abstract and represented elements of the algebra. To avoid cancellations of the denominators, the restrictions on the parameters are $\alpha+\beta\notin\mathbb{Z}$ and $\lambda_i+s+1\neq0,-1,-2,\dots $.

As in the previous subsection, to pass from the map $\Delta_{\lambda_1,\lambda_2}$ to $\Delta$, a~generalized coproduct $\mathcal{A}\to \mathcal{A}\otimes\mathcal{A}$, we express $\lambda_1$, $\lambda_2$ in terms of the central element $X$ satisfying $X(X-2)=C$ where $C$ is the Casimir element defined by~\eqref{eq:Cassl2}. This generalized coproduct is well defined if one considers the localization where $D=H^{(1)}-H^{(2)}-X^{(1)}+X^{(2)}+2\alpha+2\beta+2$ and~${H^{(i)}+X^{(i)}+2s+2}$ are invertible.

Taking the limit $\beta=+\infty$, and renaming $\gamma=\alpha+\beta+1$, reproduces the results obtained with the dual Hahn polynomials.

\section[q-Askey scheme]{$\boldsymbol{q}$-Askey scheme}\label{sec:qaskey}

This section contains the algebras and their generalized coproduct with their CG coefficients given by the finite families of polynomials of the $q$-Askey scheme.

\subsection[q-Hahn polynomials and a~deformation of the oscillator algebra]{$\boldsymbol{q}$-Hahn polynomials and a~deformation of the oscillator algebra}

Let us consider for the functions appearing in~\eqref{CG-relation} the polynomials defined by
\[Q_n(k,N;q)=
\left[\begin{array}{c} N \\ n \end{array}\right]_q
{}_3\phi_2\left(\begin{array}{c} q^{-n}, \alpha\beta q^{n-N+1}, q^{-k} \\ \alpha q, q^{-N}\end{array}\Bigg\vert q;q\right),\]
which are proportional to the standard $q$-Hahn polynomials.
Their contiguity relations called (qHIV) in~\cite{contiguity} are\footnote{We have replaced $\beta$ by $\beta q^{-N}$ compared to~\cite{contiguity}.}
\begin{gather*}
Q_n(k,N+1;q)=\frac{1-\alpha\beta q^{n+1}}{1-\alpha\beta q^{2n-N}}Q_{n-1}(k,N;q)+\frac{q^n\bigl(1-\alpha\beta q^{n-N}\bigr)}{1-\alpha\beta q^{2n-N}}Q_n(k,N;q),\\
\bigl(1-q^{N-k}\bigr)Q_n(k,N-1;q)=\frac{\bigl(1-\alpha q^{n+1}\bigr)\bigl(1-q^{n+1}\bigr)}{1-\alpha\beta q^{2n+2-N}}Q_{n+1}(k,N;q)\\
\phantom{\bigl(1-q^{N-k}\bigr)Q_n(k,N-1;q)=}{}+\frac{\alpha q^{n+1}\bigl(1-\beta q^{n+1-N}\bigr)\bigl(1-q^{N-n}\bigr)}{1-\alpha\beta q^{2n+2-N}}Q_n(k,N;q).
\end{gather*}

{\bf $\boldsymbol{q}$-oscillator algebra.} Combining the left-hand side of the second contiguity relation with condition~\eqref{eq:phimu} implies that
$\phi(\lambda,k)=1-q^k$.
Therefore, the representations of $\mathcal{A}$ are
\[H|\lambda,n\rangle=(\lambda+2n)|\lambda,n\rangle, \qquad
E|\lambda,n\rangle=|\lambda,n+1\rangle, \qquad
F|\lambda,n\rangle=(1-q^n)|\lambda,n-1\rangle.\]
Take $ K=q^{H/2}$.
In terms of this generator, the defining relations of the algebra $\mathcal{A}$ read
$
KE=q EK$, $ q KF=FK$, $ qEF-FE=q-1$.
Therefore, the algebra $\mathcal{A}$ is the $q$-oscillator algebra $osc_q$.
The Casimir element of this algebra is~${C=(1-EF)K^{-1}}$, with $C|\lambda,n\rangle=q^{-\lambda/2}|\lambda,n\rangle $.

{\bf Generalized coproduct.} As explained in the previous section, in order to determine the homomorphism $\Delta$, we look for an element $x$ of the algebra such that
\[
 x|\lambda_1,n\rangle\otimes |\lambda_2,N+1-n\rangle=\frac{1-\alpha\beta q^{n+1}}{1-\alpha\beta q^{2n-N}}|\lambda_1,n\rangle\otimes |\lambda_2,N+1-n\rangle.
\]
Knowing the actions of $K^{(1)}=K\otimes \mathbb{I}$, $K^{(2)}=\mathbb{I} \otimes K$, $C^{(1)}=C\otimes \mathbb{I}$ and $C^{(2)}=\mathbb{I} \otimes C$ on this vector, we can show that the suitable $x$ is given by
\smash{$
 x=\frac{1}{D}\bigl(1-q\alpha\beta C^{(1)}K^{(1)}\bigr)$},
where $D=1-q\alpha\beta C^{(1)}K^{(1)}\bigl(C^{(2)}K^{(2)}\bigr)^{-1}$.
A similar computation for $y$ leads to the following formula for the generalized coproduct:
\[
 \Delta(E)=\frac{1}{D}\bigl(\bigl(1-q\alpha\beta C^{(1)}K^{(1)}\bigr)(E\otimes \mathbb{I})+C^{(1)}K^{(1)}\bigl(1-q\alpha\beta \bigl(C^{(2)}K^{(2)}\bigr)^{-1}\bigr)(\mathbb{I}\otimes E)\bigr).
\]
Similarly, we compute
\[
\Delta(F)=\frac{1}{D}\bigl(\bigl(1-q\alpha C^{(1)}K^{(1)}\bigr)(F\otimes \mathbb{I})+q\alpha C^{(1)}K^{(1)} \bigl(1-\beta \bigl(C^{(2)}K^{(2)}\bigr)^{-1}\bigr)(\mathbb{I}\otimes F)\bigr).\]
This generalized coproduct is well defined in the localization of $\mathcal{A}\otimes \mathcal{A}$, where $D$
is invertible.

Keeping the same algebra $\mathcal{A}=osc_q$, generalized coproducts associated to some other polynomials of the $q$-Askey scheme can be obtained by different limits of the previous formulas since, as demonstrated in~\cite{contiguity}, the necessary contiguity relations are obtained as limits. The polynomials which may be obtained are the $q$-Krawtchouk, the quantum $q$-Krawtchouk and the affine $q$-Krawtchouk polynomials. Details are left to the interested reader.

{\bf Comparison with previous results.} In~\cite{KMM92}, a~generalized coproduct structure was given for a~deformation of the centrally extended oscillator algebra defined by~\eqref{eq:oscext}. They proved that the CG coefficients are expressed in terms of the quantum $q$-Krawtchouk. Let us mention that they also need some localizations in their algebra to be able to define their generalized coproduct which is also not coassociative. In the same spirit as in Section~\ref{sec:Hahn}, our approach allows us to obtain more general CG coefficients and we do not need to introduce an additional central charge.

\subsection[q-Racah polynomials and a~generalized coproduct for U\_q(sl\_2)]{$\boldsymbol{q}$-Racah polynomials and a~generalized coproduct for $\boldsymbol{{\rm U}_q(\mathfrak{sl}_2)}$} \label{eq:coqr}

Let us consider for the functions appearing in~\eqref{CG-relation} the polynomials defined by\footnote{We have replaced $\beta$ by $\beta q^{-N}$ and $\gamma$ by $\gamma q^N$ compared to~\cite{contiguity}.}
\[R^{(q{\rm R})}_n(k,N;q)=\left[\begin{array}{c} N \\ n \end{array}\right]_q
{}_4\phi_3\left(\begin{array}{c} q^{-n}, \alpha\beta q^{n-N+1}, q^{-k}, \gamma q^{k} \\ \alpha q, \beta\gamma q, q^{-N}\end{array}\Bigg\vert q;q\right),\]
which are proportional to the $q$-Racah polynomials.

The contiguity relations (see (RI) in~\cite{contiguity}) are
\[
R^{(q{\rm R})}_n(k,N+1;q)=\frac{1-\alpha\beta q^{n+1}}{1-\alpha\beta q^{2n-N}}R^{(q{\rm R})}_{n-1}(k,N;q)+\frac{q^n\bigl(1-\alpha\beta q^{n-N}\bigr)}{1-\alpha\beta q^{2n-N}}R^{(q{\rm R})}_n(k,N;q),
\]
and
\begin{gather*}
\bigl(1-q^{N-k}\bigr)\bigl(1-\gamma q^{N+k}\bigr)R^{(q{\rm R})}_n(k,N-1;q)\\
\qquad=\frac{\bigl(1-\beta\gamma q^{n+1}\bigr)\bigl(1-\alpha q^{n+1}\bigr)\bigl(1-q^{n+1}\bigr)}{1-\alpha\beta q^{2n+2-N}}R^{(q{\rm R})}_{n+1}(k,N;q)\\
\phantom{\qquad=}{}+\frac{\bigl(1-\beta q^{n+1-N}\bigr)\bigl(\alpha q^{n+1}-\gamma q^N\bigr)\bigl(1-q^{N-n}\bigr)}{1-\alpha\beta q^{2n+2-N}}R^{(q{\rm R})}_n(k,N;q).
\end{gather*}
Combining the left-hand side of the second contiguity relation with relation~\eqref{eq:phimu} implies that
\[\phi(\lambda_1+\lambda_2+2k,N-k)=\bigl(1-q^{N-k}\bigr)\bigl(1-\gamma q^{N+k}\bigr),\]
which implies that
$\phi(\lambda,k)=\bigl(1-q^k\bigr)\bigl(1-s q^{k+\lambda}\bigr)$,
for a~free parameter $s$ and with the constraint~${\gamma=q^{\lambda_1+\lambda_2-1}}$.
Therefore, the representations of~$\mathcal{A}$ are
\begin{gather*}
H|\lambda,n\rangle=(\lambda+2n)|\lambda,n\rangle, \qquad
E|\lambda,n\rangle=|\lambda,n+1\rangle, \\
F|\lambda,n\rangle=(1-q^n)\bigl(1- sq^{\lambda+n}\bigr)|\lambda,n-1\rangle.\end{gather*}
As previously, let us define the following generator $ K=q^{H/2}$.
In terms of this generator, the defining relations of the algebra $\mathcal{A}$ is
$
KE=qEK$, $ qKF=FK$, $ qEF-FE=(q-1)\big(1-sqK^2\big)$.
Defining \smash{$\widetilde{F}=\frac{-q^{1/2}}{(q-1)^2}K^{-1}F$} and setting $s=q^{-1}$, the above defining relations become
\[
 KE=qEK, \qquad qK\widetilde{F}=\widetilde{F}K, \qquad E\widetilde{F}-\widetilde{F}E=\frac{K-K^{-1}}{q^{1/2}-q^{-1/2}}.
\]
The previous relations are the defining relations of ${\rm U}_q(\mathfrak{sl}_2)$.
The Casimir element is
\[
 C=EFK^{-1}-q^{-1}K-K^{-1} \qquad \text{with} \quad C|\lambda,n\rangle=-\bigl(q^{\lambda/2-1}+q^{-\lambda/2}\bigr)|\lambda,n\rangle.
\]
The resulting homomorphism $\Delta_{\lambda_1,\lambda_2}\colon \mathcal{A} \to \operatorname{End}(V_{\lambda_1}\otimes V_{\lambda_2})$ is given by
\begin{gather*}
 \Delta_{\lambda_1,\lambda_2}\colon\ E\mapsto \frac{1}{D}\bigl(\bigl(1-\alpha\beta q^{1-\lambda_1/2}K^{(1)}\bigr)(E\otimes \mathbb{I})\\
 \phantom{ \Delta_{\lambda_1,\lambda_2}\colon\ E\mapsto }{}+q^{-\lambda_1/2}K^{(1)}\bigl(1-\alpha\beta q^{1+\lambda_2/2} \bigl(K^{(2)}\bigr)^{-1}\bigr)(\mathbb{I}\otimes E)\bigr),\\
 F\mapsto\frac{1}{D}\left(\frac{\bigl(1-q^{1-\lambda_1/2}\alpha K^{(1)}\bigr)\bigl(1-q^{1-\lambda_1/2}\beta\gamma K^{(1)}\bigr)}{1-sq^{1+\lambda_1/2} K^{(1)}}(F\otimes \mathbb{I})\right.\\
\left.\phantom{F\mapsto}{}+q^{1-\lambda_1/2}\alpha K^{(1)}
 \frac{\bigl(1-q^{\lambda_2/2}\beta \bigl(K^{(2)}\bigr)^{-1}\bigr)\bigl(1- q^{-\lambda_2/2}\gamma K^{(2)}/\alpha\bigr)}{1-sq^{1+\lambda_2/2} K^{(2)}}(\mathbb{I}\otimes F)\right)
\end{gather*}
with $D=1-q^{1-\lambda_1/2+\lambda_2/2}\alpha\beta K^{(1)}\bigl(K^{(2)}\bigr)^{-1}$.
The same notations are kept for the abstract and represented elements of the algebra. As for the Racah case (see Section~\ref{ssec:Racah}), we can pass from the previous map $\Delta_{\lambda_1,\lambda_2}$ to $\Delta$, a~generalized coproduct $\mathcal{A}\to \mathcal{A}\otimes\mathcal{A}$, well defined only in a~localization.

Taking $s=q^{-1}$, $\beta=0$, $\alpha=q^{\lambda_1-1}$ and $\gamma=q^{\lambda_1+\lambda_2-1}$, the CG coefficients are associated to the dual $q$-Hahn polynomials and $\Delta_{\lambda_1,\lambda_2}$ simplifies to
\begin{gather*}
 \Delta_{\lambda_1,\lambda_2}\colon\ E\mapsto E\otimes \mathbb{I} +q^{-\lambda_1/2} K\otimes E,\qquad
  F\mapsto F\otimes \mathbb{I} +q^{\lambda_1/2} K\otimes F.
\end{gather*}
With the Drinfeld twist $\mathscr{F}=\mathbb{I} \otimes K^{\lambda_1/2}$, the twisted coproduct $\widetilde{\Delta}_{\lambda_1,\lambda_2}=\mathscr{F}\Delta_{\lambda_1,\lambda_2} \mathscr{F}^{-1}$ is
\begin{gather*}
 \widetilde{\Delta}_{\lambda_1,\lambda_2}\colon\ E\mapsto E\otimes \mathbb{I} + K\otimes E,\qquad
  F\mapsto F\otimes \mathbb{I} + K\otimes F.
\end{gather*}
One recognizes the usual coproduct of ${\rm U}_q(\mathfrak{sl}_2)$. We recover that the dual $q$-Hahn polynomials are the CG coefficients of ${\rm U}_q(\mathfrak{sl}_2)$ with the usual coproduct~\cite{Arv}.

Keeping the same algebra $\mathcal{A}={\rm U}_q(\mathfrak{sl}_2)$, generalized coproducts associated to other polynomials of the $q$-Askey scheme (dual $q$-Hahn and dual $q$-Krawtchouk) can be obtained by different limits of the previous formulas since, as demonstrated in~\cite{contiguity}, the necessary contiguity relations are obtained as limits. Details are left to the interested reader.

\section{Outlook}\label{sec:outlook}

In this paper, we have introduced several new generalized coproducts on well-known algebras.
A natural line of inquiry is whether these generalized coproducts are related to existing ones, for instance through a~Drinfeld twist. Specifically, for the one presented in Section~\ref{eq:coqr} (which extends the known result for ${\rm U}_q(\mathfrak{sl}_2)$), can it be proven that it is obtainable from the usual ${\rm U}_q(\mathfrak{sl}_2)$ coproduct via a~Drinfeld twist? In~\cite{GZ92}, the $q$-Racah polynomials appear as ``twisted'' Clebsch--Gordan coefficients for ${\rm U}_q(\mathfrak{sl}_2)$, i.e., they are the overlap coefficients between rotated bases. It would be interesting to investigate the relationship between the rotation introduced in~\cite{GZ92} and the generalized coproduct defined in this paper. Specifically, we believe that this rotation can be interpreted as a~Drinfeld twist of our map $\Delta$, recovering the standard coproduct of ${\rm U}_q(\mathfrak{sl}_2)$.
Similarly, the twisted primitive elements of ${\rm U}_q(\mathfrak{sl}_2)$, which have been introduced in~\cite{Koelink96, Koor93,NM90} to study zonal spherical functions on the ${\rm SU}_q(2)$ quantum group, could be given an interpretation of the generalized coproduct introduced here. It is also known that $q$-Racah polynomials appear naturally in the context of the dynamical quantum Yang--Baxter equation~\cite{KR01} and of the reflection equation~\cite{Bas}. Consequently, it is natural to ask whether the map $\Delta$ obtained in this work can also be derived within those frameworks.

Furthermore, since the generalized coproducts are generally non-cocommutative, it would be valuable to investigate whether an associated universal $R$-matrix exists. A next step is to examine the coassociativity of these generalized coproducts. If they are not coassociative, computing the associated coassociator would be a~worthwhile endeavor, see, for example, Remark~\ref{rema_coassociativitypq}.

The CG coefficients are proportional to the $3j$-symbols. For $\mathfrak{sl}_2$, it is possible to compute the $3nj$-symbols, which govern the decomposition of tensor products involving more than two irreducible representations into a~direct sum of irreducible representations. These higher-order symbols are also known to be associated with various special functions. We anticipate that similar computations will be possible for the homomorphisms $\Delta$ introduced here, leading to interesting special functions.
Finally, we focus on $\mathfrak{osc}$, $\mathfrak{sl}_2$ or ${\rm U}_q(\mathfrak{sl}_2)$, but we believe that our approach can also be useful to obtain generalized coproduct for other Lie algebras as $\mathfrak{su}(1,1)$ generalizing previous results obtained in this context~\cite{Alvarez}.
We plan to address these open problems in future work.

\subsection*{Acknowledgments} The authors want to thank warmly the referees for their valuable comments.
NC is partially supported by the international research project AAPT of the CNRS.
 LV is funded in part through a~discovery grant of the Natural Sciences and Engineering Research Council (NSERC) of Canada.

\pdfbookmark[1]{References}{ref}
\LastPageEnding

\end{document}